\documentclass{fic-l}
\usepackage{amsmath}
\usepackage{amsfonts}
\usepackage{amssymb}
\usepackage{amscd}
\usepackage{amsthm}
\usepackage{latexsym}
\newtheorem{theorem}{Theorem}[section]
\newtheorem{lemma}[theorem]{Lemma}
\newtheorem{proposition}[theorem]{Proposition}
\newtheorem{corollary}[theorem]{Corollary}

\theoremstyle{definition}

\theoremstyle{remark}

\numberwithin{equation}{section}

\newcommand{\bt}{\blacktriangleright}
\newcommand{\ot}{\otimes}

\newcommand{\ov}{\overline}
\newcommand{\h}{\mathcal{H}}
\newcommand{\sw}[1]{^{(#1)}}

\begin{document}
\title
{On the Cyclic Homology of Hopf Crossed Products}

\author {M. Khalkhali}
\address{Department of Mathematics \\ University of Western Ontario\\
London, ON, Canada} 
\author{ B. Rangipour}
\address{Department of Mathematics \\ University of Western Ontario\\ London, ON, Canada} 
\subjclass{Primary 16W30; Secondary 18G60, 19D55}

\begin{abstract}
We consider  Hopf crossed products of the the type $A\#_\sigma \mathcal{H}$, where $\mathcal{H}$ is a cocommutative 
Hopf algebra, $A$ is an $\mathcal{H}$-module algebra and $\sigma$ is a ``numerical" convolution  invertible 2-cocycle 
on $\mathcal{H}$. we give an spectral sequence that converges to the cyclic homology of $A\#_\sigma \mathcal{H}$ and 
 identify the $E^1$ and $E^2$ terms of the spectral sequence. 
\end{abstract}
\maketitle
\section{Introduction}
A celebrated problem in noncommutative geometry, more precisely in cyclic homology theory,
 is to compute the cyclic homology of a crossed product algebra. The interest in this problem stems  from 
  the fact that, according to a guiding principle in noncommutative geometry \cite{ac}, crossed products play the role 
  of    ``noncommutative quotients" in situations where the  usual set theoretic quotients are ill behaved. 
   For  example, when a (locally compact) group $G$ acts on a locally compact 
   Hausdorff space $X$, the quotient space $X/G$  
   may not be well behaved, e.g. may not be even a Hausdorff space. 
   The crossed product $(C^*-)$ algebra $C_0(X)\ltimes G$, however, is a good replacement
    for $X/G$ \cite{ac}. In fact, if the action of $G$ is
     free and proper, then  by a theorem of Rieffel \cite{r} the $C^\ast$-algebra 
      of continuous functions vanishing at infinity on $X/G$, denoted by $C_0(X/G)$ is in a 
     suitable $C^*$-algebraic sense, Morita equivalent to the crossed product algebra $C_0(X)\ltimes G$. Since 
      $K$-Theory, Hochschild homology and  cyclic homology are Morita invariant functors, 
      replacing the commutative algebra $C_0(X/G)$ by the 
      noncommutative algebra $C_0(X)\ltimes G$ results in no loss of information.
      
       It is 
      therefore natural and desirable to develop tools to compute the cyclic 
      homology of crossed product algebras. Most of the results obtained 
      so far are concerned  with the  action of groups  on algebras \cite{ft,gj,n}. 
For  Hopf algebra crossed products,  \cite{gugu}  gives a complete answer for 
Hochschild homology but it is not clear how to extend its method to  cyclic homology.
 If one is only interested in smash products as
 opposed to crossed products, one can find a complete answer in \cite{ak} in terms of a spectral
sequence converging to the cyclic homology of the smash product algebra $A\# \mathcal{H}$.
     
The goal of this article  is to extend the results of \cite{ak} to Hopf algebra 
crossed products. Due to the fact that very complicated formulas appear  in our 
constructions, however, we had to assume 
that the Hopf algebra is  cocommutative and the  2-cocycle takes values in the ground field. 
Under these conditions,  we give  a spectral 
sequence for the  cyclic homology of a crossed product  algebra $A\underset{\sigma}{\#} \mathcal{H}$,
when the cocycle $\sigma$ is convolution  invertible and takes values in the ground  field $k$. 
The method of proof 
is similar to the  one used in \cite{ak} which is based on the generalized cyclic Eilenberg-Zilber  theorem of 
\cite{gj}.  Though we think the same method should apply to arbitrary 
 crossed products $A\underset{\sigma}{\#} \mathcal{H}$  (with convolution 
 invertible cocycles), due to technical  difficulties  we are not able to verify this. 
 
 One of our main motivations to consider cocycle crossed products is to find simple 
  methods to compute the cyclic cohomology of  `` noncommutative toroidal orbifolds" considered in \cite[Sec. 9]{s} and \cite{d}.
    These 
  examples are suggested by applications of noncommutative geometry to string theory and M(atrix) theory. 
  In these examples one considers algebras of the type
   $B^d_{\theta,\sigma}=A^d_\theta\underset{\sigma}{\#}\mathbb{C}G$, 
   where $G$ is a finite group acting by automorphisms on the noncommutative $d$-dimensional torus $A^d_\theta$ and 
   $\sigma\in H^2(G,U(1))$ is a group $2$-cocycle on $G$. Corollary \ref{cor} shows that the cyclic cohomology 
    of $B^d_{\theta,\sigma}$  can always be computed from a cyclic complex much simpler than  the original 
   cyclic complex of the algebra $B^d_{\theta,\sigma}$. 
\section{Preliminaries }\label{sec2}
In this paper we work over a fixed ground field $k$. All algebras are unital  
associative algebras over $k$  
and all
modules are unitary. The unadorned tensor product $\ot$ means tensor product  over $k$.
We denote the coproduct of a Hopf algebra by $\Delta$, the counit by $\epsilon$
 and the antipode by $S$.
We use Sweedler's notation and write  $\Delta(h)=h^{(1)}\ot h^{(2)}$ to denote 
the coproduct, where summation is understood. Similarly,
 we write $\Delta^{(n)}(h)=h^{(1)}\ot h^{(2)}\ot \dots \ot h^{(n+1)}$ to denote the iterated 
 coproducts defined by  $\Delta^{(1)}=\Delta$ and $\Delta^{(n+1)}=(\Delta\ot 1)\circ \Delta^{(n)}.$
 
 We recall the concept of Hopf crossed product,  introduced for the first time in \cite{dt}, 
 and independently \cite{mo}. A good reference for this notion is Chapter 7 of \cite{m}.   Let $\mathcal{H}$ be a Hopf algebra and $A$ an algebra. 
Recall from \cite{mo} and \cite{dt} that a {\it weak action}  of $\mathcal{H}$  on $A$ is a linear 
map $\mathcal{H}\ot A \longrightarrow A$, $h\ot a\rightarrow h(a)$ such that,
for all  $h\in\mathcal{H}$, and  $a,b\in A$
\begin{itemize}

\item[1)] ~ $h(ab)=h^{(1)}(a)h^{(2)}(b)$,
\item[2)]~  $h(1)=\epsilon(h)1$,
\item[3)]~  $1(a)=a$.
\end{itemize}

    By an {\it action } of $\mathcal{H}$ on $A$, we  mean a weak  action 
    such that $A$ is an $\mathcal{H}$-module, i.e. for all $h,l\in \mathcal{H}$ and $a\in A$ we
     have $h(l(a))=hl( a)$. In the latter  case we say $A$ is an $\mathcal{H}-${\it module algebra}.
     
     Let $A$ be an $\mathcal{H}$-module algebra. The {\it smash product} $A\#\mathcal{H}$ of $A$
and $\mathcal{H}$ is an associative algebra whose underlying vector space is $A\ot \mathcal{H}$ and
 whose multiplication is defined by 
 $$(a\ot h)(b\ot l)= ah^{(1)}( b) \ot h^{(2)}l.$$
 
 If, on the other hand, we have only a weak  action of $\mathcal{H}$ on $A$ the above formula does
  not define an associative multiplication, and a modification is needed. Given a linear map
  $\sigma:\mathcal{H}\ot \mathcal{H}\rightarrow A$ one defines a (not necessarily unital or
  associative) multiplication on $A\ot\mathcal{H}$  by \cite{mo,dt} 
  $$(a\ot h)(b\ot l)=ah^{(1)}(b)\sigma(h^{(2)},l^{(1)})\ot h^{(3)}l^{(2)}.$$ 
   It can be shown that the above formula defines an associative product with $1\ot 1$  as its unit, 
   if and only if $\sigma$
    and the weak action enjoy  the following properties:
    \begin{itemize}

\item[1)]~(Normality) For all $h\in\mathcal{H}$,  $\sigma(h,1)=\sigma(1,h)=\epsilon(h)1$.
\item[2)]~(Cocycle property) For all $h,l,m\in \mathcal{H}$, \\
$                  \sum h^{(1)}(\sigma(l^{(1)}
,m^{(1)}))\sigma(h^{(2)},l^{(2)}m^{(2)})=\sum \sigma(h^{(1)},l^{(1)})\sigma(h^{(2)}l^{(2)},m),$
\item[3)](Twisted module property) For all $h,l\in \mathcal{H}$ and $a\in A$, \\
$\sum h^{(1)}(l^{(1)}(a))\sigma(h^{(2)},l^{(2)})=\sum \sigma(h^{(1)},l^{(1)})l^{(2)}h^{(2)}(a)).$
\end{itemize}

The cocycle $\sigma$ is said to be  {\it convolution invertible } if it is an invertible element 
of the convolution algebra $Hom_k(\mathcal{H}\ot \mathcal{H},A)$. 
Now assume the Hopf algebra $\mathcal{H}$
 is  cocommutative, $\sigma: \mathcal{H}\ot \mathcal{H}\longrightarrow k1_A$ 
 takes values in the ground field 
 $k$, and $\sigma$ is invertible. Then it follows that A is an $\mathcal{H}$-module algebra, i.e. the 
  weak action in the above  situation is in fact an   action. 
 To prove this let $a\in A$ be fixed. Define two functions in
  $Hom_k(\mathcal{H}\ot \mathcal{H}, A)$ by 
 $F(h,l)=\sum h^{(1)}(l^{(1)}(a))\sigma(h^{(2)},l^{(2)})$ and 
 $G(h,l)=\sum \sigma(h^{(1)},l^{(1)})l^{(2)}h^{(2)}(a).$ 
  Then $F=G$ by  the twisted module property of $\sigma$, so  $F\ast \sigma^{-1}=G\ast\sigma^{-1} $.
 In other word  
     $$h(l(a))=F\ast\sigma^{-1}(h,l)=G\ast\sigma^{-1}(h,l)=hl(a).$$
	 One notes that the above proof remains valid when $\mathcal{H}$ is cocommutative and $\sigma$ takes its 
	 values in the center of $A$, instead of $k$.  
   
One of the main technical tools used   in \cite{ak}  to derive a spectral
 sequence for the cyclic homology of smash 
 products is the generalized cyclic
 Eilenberg-Zilber theorem. This result was  first stated in \cite{gj} but 
 its first algebraic  proof appeared in \cite{kr1}. The idea 
 of using  an Eilenberg-Zilber type theorem to derive  a spectral sequence for 
 cyclic homology of smash products (for the action of groups) is due to Getzler and Jones \cite{gj}.
  We find it remarkable that the same idea works in the case of Hopf algebra 
  crossed product (with convolution invertible cocycle). 
  In the following we recall the definitions of (para)cyclic modules, cylindrical modules
   and  state the Eilenberg-Zilber theorem for cylindrical modules.
   
   Recall that a {\it paracyclic module } is a simplicial $k$-module $\{M_n\}_{n\ge 0}$,
    such that the following extra relations are satisfied \cite{ft,gj}:

    \begin{align*}
     & \delta_i\tau=\tau\delta_{i-1}, &&\delta_0\tau=\delta_n, &&& 1\le i\le n,\\
     &\sigma_i\tau=\tau\sigma_{i+1}, &&\sigma_0\tau=\tau^2\sigma_n, &&& 1\le i\le n,
     \end{align*}
     where $\delta_i:M_n\rightarrow M_{n-1}$, ~$\sigma_i:M_n\rightarrow M_{n+1}$, $0\le i\le n$, 
     are faces and degeneracies of the simplicial module $\{M_n\}_{n\ge 0}$ and $\tau: M_n\rightarrow M_n$, 
     $n\ge 0$ are $k$-linear maps. 
     If furthermore we have $\tau^{n+1}=id_{M_n}$,
     for all $n\ge 0$, then we say that we have a {\it cyclic module}. 
     
     We denote the cyclic module of an associative unital $k$-algebra $A$, by $A^\natural$. It is defined 
      by $A^\natural_n=A^{\ot(n+1)}$, $n\ge 0$, and  simplicial and cyclic operations defined by       
   \begin{eqnarray*}
\delta_i(a_0\otimes\dots \otimes a_n)&=&a_0\otimes a_1\otimes \dots \otimes a_ia_{i+1}  \dots \otimes a_n,\;\;\;0\leq
 i\leq n-1,\\
\delta_{n}(a_0\otimes\dots \otimes a_n)&=&a_na_0\otimes a_1\otimes \dots \otimes a_{n-1},\\
\sigma_i(a_0\otimes \dots\otimes a_n)&=& a_0\otimes \dots\otimes a_i\otimes 1\dots   \otimes a_{n},\;\;\;0\leq i\leq
 n,\\
\tau_n(a_0\otimes \dots\otimes a_n)&=& a_n\otimes a_0 \dots\otimes a_{n-1}.
\end{eqnarray*} 
     
To any cyclic module one associates its cyclic homology groups \cite{ld,ac88}.
In particular, the cyclic homology groups of $A^\natural$ are denoted by  $HC_n(A)$, $n\ge 0$,
and  are called cyclic homology of $A$.  
     
     By a {\it biparacyclic} module we mean a doubly graded sequence of $k$-modules $\{M_{p,q}\}_{p,q\ge 0}$     
      such that each row and each column is  a paracyclic module and all vertical 
      operators commute with all horizontal operators. 
      In particular a {\it bicyclic } module is a biparacyclic module such that each row and 
      each column is a cyclic module. We denote the horizontal and vertical 
      operators of a biparacyclic module by $(\delta_i,\sigma_i,\tau)$ and $(d_i,s_i,t)$ respectively.      
      By a {\it cylindrical module } we mean a biparacyclic module such that for all $p,q\ge 0$, 
      \begin{equation}
      \tau^{p+1}t^{q+1}=id_{M_{p,q}}.\label{*}      
      \end{equation}
Given a cylindrical module $M$, its {\it diagonal}, denoted by $dM$, is a 
cyclic module defined by $(dM)_n=M_{n,n}$ and with simplicial and cyclic operators
given by $\delta_id_i$, $\sigma_is_i$, and $\tau t$.  In view of (\ref{*}), it  is a cyclic module.
 The total complex of a cylindrical module, denoted by $Tot(M)$,
is a mixed complex with operators given by $b+\bar{b}$ and $B+T\bar{B}$, where $T=1-(bB+Bb)$.
 Here $b$ (resp. $\bar{b}$) and $B$ (resp. $\bar{B}$) are the vertical (resp. horizontal) Hochschild and Connes
  boundary operators of cyclic modules.  
Note that it differs
from the usual notion  of total complex in that we use $B+T\bar{B}$ instead of $B+\bar{B}$. In 
fact the latter choice won't give us a  mixed complex \cite{gj}. It can be checked that $Tot(M)$ is a mixed complex.
 Given a cylindrical or cyclic module $M$,  we denote its normalization by $N(M)$. 
          
 The following theorem is the main technical  result that enables  us to derive spectral 
 sequences for the cyclic homology of crossed product algebras.
\begin{theorem}[Generalized cyclic Eilenberg-Zilber theorem (\cite{kr1,gj})]
For any  cylindrical module $M$ 
 there is a natural quasi-isomorphism of mixed complexes 
$f_0+uf_1:Tot(N(M))\longrightarrow N(dM)$, where $f_0$
is the shuffle map.
           
\end{theorem}

\section{ A Spectral sequence for Hopf crossed products}
  Let $\mathcal{H}$ be a cocommutative Hopf algebra, $A$ a  
  left $\mathcal{H}$- module algebra  and
   $\sigma: \mathcal{H}\ot\mathcal{H}\longrightarrow k1_A $ a two cocycle satisfying 
   the cocycle  conditions  1), 2), and 3) in Section \ref{sec2}. We further assume that 
   $\sigma$ is convolution invertible.  
    We introduce a cylindrical module 
    $$A\natural_\sigma\mathcal{H}=\{\mathcal{H}^{\ot(p+1)}\ot A^{\ot(q+1)}\}_{p,q\ge 0}$$
     with  vertical and horizontal  simplicial and cyclic operators 
     $(\delta,\sigma,\tau)$ and $(d,s,t)$, defined as follows

\begin{multline*}
\tau( g_0 , \dots , g_p \mid a_0 , \dots , a_q) 
= (g_0^{(2)}, \dots , g_p^{(2)} \mid \\
\shoveright{ S(g_0^{(1)}\dots g_p^{(1)})(a_q), a_0,a_1, \dots ,a_{q-1})} \\
\shoveleft{
\delta_i(g_0,\dots,g_p \mid a_0 , \dots , a_q)= (g_0,\dots,g_p \mid a_0,\dots , a_i a_{i+1},\dots , a_q) 
\;\;\; 0 \le i <q } \\
\shoveleft{\delta^{p,q}_q (g_0,\dots,g_p \mid a_0 , \dots , a_q)= (g_0^{(2)},\dots,g_p^{(2)}
 \mid }\\ 
 \shoveright{ S(g_0^{(1)}\dots g_p^{(2)})(a_q)a_0
 ,a_1\dots , a_{q-1})} \\
\shoveleft{\sigma_i(g_0,\dots,g_p \mid a_0 , \dots , a_q) = (g_0,\dots,g_p 
\mid a_0,\dots , a_i , 1 ,
 a_{i+1},\dots , a_q) 
\;\; 0 \le i \le q }\\
\shoveleft{t( g_0 , \dots , g_p \mid a_0 , \dots , a_q)
=(g_p^{(q+2)},g_0, \dots , g_{p-1} \mid g_p^{(1)} (a_0), \dots
,g_p^{(q+1)}( a_{q}))} \\
\shoveleft{d_i(g_0,\dots,g_p \mid a_0 , \dots , a_q)
=(g_0,\dots,g_i^{(1)} g_{i+1}^{(1)},\dots,g_p \mid} \\
\shoveright{\sigma(g_i^{(2)},g_{i+1}^{(2)})a_0,
\dots , a_q) \;\;\; 0 \le i <p  }\\ 
\shoveleft{d_p (g_0,\dots,g_p \mid a_0 , \dots , a_q)
= (g_p^{(q+2)}g_0^{(1)},g_1,\dots,g_{p-1} \mid }\\
\shoveright{\sigma(g_p^{(q+3)},g_0^{(2)}) g_p^{(1)} 
(a_0),\dots , g_p^{(q+1)}( a_{q}))}\\
\shoveleft{s_i (g_0,\dots,g_p \mid a_0 , \dots , a_q)
=(g_0,\dots,g_i,1,g_{i+1},\dots,g_p \mid a_0,\dots , a_q) \;\;\; 0 \le i \le p. } 
\end{multline*} 
\begin{theorem}
Endowed with the above operators, $A\natural_\sigma\mathcal{H}$ is a cylindrical module.
\end{theorem}
\begin{proof}
 We should check that every  row and every  column is a paracyclic module, 
  and   vertical operator commutes  with each   horizontal operator, and  in addition the identity (\ref{*}) holds.
 Since the weak action in our situation is actually an action and the vertical  operators are the same as
 the vertical  operators in (\cite{ak} Theorem 3.1 ), we refer the reader to \cite{ak} 
 for the proof  that the columns  form paracyclic modules. 
 To check that the rows  are paracyclic modules we need to verify the following identities

\begin{align*}
d_id_j &= d_{j-1}d_i \hspace{15pt} i<j\\
s_is_j &= s_{j+1}s_i \hspace{15pt} i\le j \\
d_is_j &=
\begin{cases}
s_{j-1}d_i &\text{$i<j$}\\
\mbox{identity}      &\text{$i=j$ or $i=j+1$}\\
s_jd_{i-1} &\text{$i>j+1$}.
\end{cases}
\end{align*}

\begin{align*}
d_it_n = t_{n-1}d_{i-1} \hspace{15pt} 1\le i \le n\; , \quad
d_0 t_n = d_n  \\
s_i t_n =  t_{n+1} s_{i-1} \hspace{15pt} 1\le i\le n\; , \quad
s_0 t_n = t_{n+1}^2 s_n 
\end{align*}

 We just check $d_id_{i+1}=d_id_i$ and the cylindrical  module condition (\ref{*}).
  The rest  can be proved by the same techniques. We have:
 \begin{align*}
 &d_{i}d_{i+1}(g_0,\dots,g_p \mid a_0 , \dots , a_q)=\\
 &d_{i}(g_0,\dots,g_{i+1}^{(1)}g_{i+2}^{(1)},\dots g_p 
 \mid\sigma(g_{i+1}^{(2)},g_{i+2}^{(2)}) a_0, \dots , a_q)=\\
 &(g_0,\dots,g_i^{(1)}g_{i+1}^{(1)}g_{i+2}^{(1)},\dots g_p\mid
  \sigma(g_i^{(2)},g_{i+1}^{(2)}g_{i+2}^{(2)})\sigma(g_{i+1}^{(3)},g_{i+2}^{(3)})a_0\dots,a_q),\\ 
  &\text{that by using the cocycle property 2) in Section 2  is equal to }\\
  &(g_0,\dots,g_i^{(1)}g_{i+1}^{(1)}g_{i+2}^{(1)},\dots g_p\mid
  \sigma(g_i^{(2)},g_{i+1}^{(2)})\sigma(g_i^{(3)}g_{i+1}^{(3)},g_{i+2}^{(3)})a_0\dots,a_q)=\\
 & d_id_i(g_0,\dots,g_p \mid a_0 , \dots , a_q).  
 \end{align*}
 
 Next we check the cylindrical module condition (\ref{*}).
 We have :\\
 $t^{p+1} \tau^{q+1}(g_0,\dots,g_p \mid a_0, \dots , a_q)$
\begin{align*}
&= t^{p+1} \tau^q(g_0^{(2)},\dots , g_p^{(2)} \mid S(g_0^{(1)} g_{1}^{(1)} 
\dots g_p^{(1)}) \cdot a_q,a_0, \dots , a_{q-1})\\ 
&= t^{p+1}(g_0^{(q+1)},\dots,g_p^{(q+1)} \mid S(g_0^{(q)} \dots g_p^{(q)}) 
\cdot a_0 ,S(g_0^{(q-1)} \dots g_p^{(q-1)}) \cdot a_1 ,\\
&  \hspace{9.5cm} \dots ,S(g_0^{(0)} \dots g_p^{(0)}) \cdot a_q )\\
&= t^{p}(g_p^{(2q+2)},g_0^{(q+1)},\dots,g_{p-1}^{(q+1)} \mid (g_p^{(q+1)} S(g_0^{(q)} \dots g_p^{(q)}) 
)\cdot a_0 ,\\
&  \hspace{3cm} (g_p^{(q+2)}S(g_0^{(q-1)} \dots g_p^{(q-1)}) \cdot a_1 ,
\dots ,(g^{(2q+1)}S(g_0^{(0)} \dots g_p^{(0)})) \cdot a_q )\\
&= t^{p}(g_p^{(2q)},g_0^{(q+1)},\dots,g_{p-1}^{(q+1)} \mid (S(g_0^{(q)} \dots g_{p-1}^{(q)}) 
)\cdot a_0 ,\\
&  \hspace{3.5cm} (g_p^{(q)}S(g_0^{(q-1)} \dots g_p^{(q-1)}) \cdot a_1 ,
\dots ,(g_p^{(2q-1)}S(g_0^{(0)} \dots g_p^{(0)})) \cdot a_q )\\
&= t^{p}(g_p ,g_0^{(q+1)},\dots,g_{p-1}^{(q+1)} \mid (S(g_0^{(q)} \dots g_{p-1}^{(q)}) 
)\cdot a_0 ,\\
&  \hspace{4.5cm} (S(g_0^{(q-1)} \dots g_{p-1}^{(q-1)}) \cdot a_1 ,
\dots ,(S(g_0^{(0)} \dots g_{p-1}^{(0)})) \cdot a_q )\\
&=(g_0,\dots,g_p \mid a_0, \dots , a_q). 
\end{align*}
The theorem is proved. 
\end{proof}

Next we show that the diagonal  of the above cylindrical module, $d(A\natural_\sigma\mathcal{H})$,
is isomorphic with  the cyclic module  $(A\#_\sigma\mathcal{H})^\natural$ 
associated with the crossed product algebra. To this end we define maps 
 $\Phi: (A\#_\sigma \mathcal{H})^\natural\rightarrow d(A\natural_\sigma\mathcal{H}) $ and 
 $\Psi: d(A\natural_\sigma\mathcal{H})\rightarrow   (A\#_\sigma \mathcal{H})^\natural$ by the following formulas\\
$\Phi ( a_0 \otimes g_0,\dots ,a_n \otimes g_n )=$
\begin{eqnarray*}
& & (g_0^{(2)},g_1^{(3)},\dots,g_n^{(n+2)} \mid  S(g_0^{(1)}g_1^{(2)}
 \dots g_n^{(n+1)}) \cdot a_0,S(g_1^{(1)}g_2^{(2)} \dots g_n^{(n)}) \cdot a_1,\dots \\
& & \hspace{5cm} ,
S(g_{n-1}^{(1)}g_n^{(2)})\cdot a_{n-1},S(g_n^{(1)}) \cdot a_n),   
\end{eqnarray*}
$\Psi (g_0, \dots , g_n \mid a_0 , \dots ,a_n)=$  
\begin{eqnarray*}
((g_0^{(1)} g_1^{(1)} \dots g_n^{(1)}) \cdot a_0 \otimes g_0^{(2)} ,(g_1^{(2)}
 \dots g_n^{(2)}) \cdot a_1 \otimes g_1^{(3)},
\dots , g_n^{(n+1)} \cdot a_n \otimes g_n^{(n+2)}).
\end{eqnarray*}
\begin{theorem}
 The above maps, $\Phi$, $\Psi$,  are morphisms of cyclic modules and are inverse  to  one another.
  \end{theorem}
  \begin{proof}
  It is not hard to see that $\Phi$ and $\Psi$ are inverse of each other. 
  We just prove that $\Phi$ is a cyclic map.
  We first  verify the commutativity of $\Phi$ and the cyclic operators, i.e., the relation  
  $ t \tau \Phi = \Phi \tau_{A \#_\sigma \mathcal{H}}.$ We have:  
\begin{align*}
& (t \tau) \Phi (a_0 \otimes g_0, \dots , a_n \otimes g_n)\\
&=t\tau(g_0^{(2)},g_1^{(3)}, \dots , g_n^{(n+2)} \mid S(g_0^{(1)} g_1^{(2)} 
\dots g_n^{(n+1)})\cdot a_0 , S(g_1^{(1)} g_2^{(2)} 
\dots g_n^{(n)})\cdot a_1,\\ 
&  \hspace{6cm} \dots,S(g_{n-1}^{(1)} g_n^{(2)}) \cdot a_{n-1} ,S(g_n^{(1)}) \cdot a_n )\\
&= t (g_0^{(3)},g_1^{(4)}, \dots , g_n^{(n+3)} 
\mid S(g_0^{(2)} g_1^{(3)} \dots g_n^{(n+2)})S(g_n^{(1)})\cdot a_n ,\\
&  \hspace{4cm} S(g_0^{(1)} g_1^{(2)} 
\dots g_n^{(n+1)})\cdot a_0 \dots,S(g_{n-1}^{(1)} g_n^{(2)}) \cdot a_{n-1}  )\\
&= (g_n^{((2n+4)}, g_0^{(3)} ,g_1^{(4)}, \dots , g_{n-1}^{(n+2)} 
\mid g_n^{(n+3)} S(g_0^{(2)} g_1^{(3)} \dots g_n^{(n+1)})S(g_n^{(1)})\cdot a_n ,\\
&  \hspace{2cm} g_n^{(n+4)} S(g_0^{(1)} g_1^{(2)} 
\dots g_n^{(n+1)})\cdot a_0 \dots,g_n^{(2n+3)} S(g_{n-1}^{(1)} g_n^{(2)}) \cdot a_{n-1}  )\\
&= (g_n^{((2n+3)}, g_0^{(3)} ,g_1^{(4)}, \dots , g_{n-1}^{(n+2)}
 \mid \epsilon (g_{n}^{(n+2)}) S(g_n^{(1)} g_0^{(2)} \dots g_{n-1}^{(n+1)})\cdot a_n ,\\
&  \hspace{2cm} g_n^{(n+3)} S(g_0^{(1)} g_1^{(2)} 
\dots g_n^{(n+1)})\cdot a_0 ,\dots,g_n^{(2n+2)} S(g_{n-1}^{(1)} g_n^{(2)}) \cdot a_{n-1}  )\\
&= (g_n^{((2n+2)}, g_0^{(3)} ,g_1^{(4)}, \dots , g_{n-1}^{(n+2)}
 \mid  S(g_n^{(1)} g_0^{(2)} \dots g_{n-1}^{(n+1)})\cdot a_n ,\\
&  \hspace{1.5cm} \epsilon (g_{n}^{(n+1)}) S(g_0^{(1)} g_1^{(2)} 
\dots g_{n-1}^{(n)})\cdot a_0 ,\dots,g_n^{(2n)} S(g_{n-1}^{(1)} g_n^{(2)}) \cdot a_{n-1}  )\\
&= (g_n^{(2)}, g_0^{(3)} ,g_1^{(4)}, \dots , g_{n-1}^{(n+2)} 
\mid  S(g_n^{(1)} g_0^{(2)} \dots g_{n-1}^{(n+1)})\cdot a_n ,\\
&  \hspace{3cm}  S(g_0^{(1)} g_1^{(2)} \dots g_{n-1}^{(n)}) \cdot a_0 , \dots 
,\dots,\epsilon(g_n^{(2)}) S(g_{n-1}^{(1)})  \cdot a_{n-1}  )\\
&= (g_n^{(2)}, g_0^{(3)} , \dots , g_{n-1}^{(n+2)} \mid  S(g_n^{(1)} g_0^{(2)} \dots g_{n-1}^{(n+1)})\cdot a_n ,\\
&  \hspace{4cm}  S(g_0^{(1)} g_1^{(2)} \dots g_{n-1}^{(n)}) \cdot a_0 , \dots 
,\dots, S(g_{n-1}^{(1)})  \cdot a_{n-1}  )\\
&= \Phi (a_n \otimes g_n,a_0 \otimes g_0 , \dots , a_{n-1} \otimes g_{n-1})
= \Phi (\tau_{A \#_\sigma \mathcal{H}}( a_0 \otimes g_0, \dots , a_n \otimes g_n)).
\end{align*}

Next  we check the commutativity of $\Phi$ and the face operators, i.e. the relation 
 $d_i\delta_i\Phi=\Phi d_i^{A\#_\sigma\mathcal{H}}$. For $0\le i<n,$ we have: 
\begin{align*}
&d_i\delta_i\Phi( a_0 \otimes g_0, \dots , a_n \otimes g_n)= \\
& d_i\delta_i(g_0^{(2)},g_1^{(3)}, \dots , g_n^{(n+2)} \mid S(g_0^{(1)} g_1^{(2)} 
\dots g_n^{(n+1)})\cdot a_0 , S(g_1^{(1)} g_2^{(2)} 
\dots g_n^{(n)})\cdot a_1,\\ 
&  \hspace{6cm} \dots,S(g_{n-1}^{(1)} g_n^{(2)}) \cdot a_{n-1} ,S(g_n^{(1)}) \cdot a_n )=\\
&d_i((g_0^{(2)},g_1^{(3)}, \dots , g_n^{(n+2)} \mid S(g_0^{(1)} g_1^{(2)} 
\dots g_n^{(n+1)})\cdot a_0 , S(g_1^{(1)} g_2^{(2)} 
\dots g_n^{(n)})\cdot a_1,\\ 
& \dots,S(g_{i+1}^{(1)}, \dots,  g_n^{(n+1-i)})
\cdot(S(g_i^{(1)})(a_i)a_{i+1}), \dots,  S(g_{n-1}^{(1)} g_n^{(2)}) \cdot a_{n-1} ,S(g_n^{(1)}) \cdot a_n ))=\\
&((g_0^{(2)},g_1^{(3)}, \dots ,g_i^{(i+2)}g_{i+1}^{(i+3)},\dots 
 g_n^{(n+2)} \mid \sigma(g_i^{(i+3)},g_{i+1}^{(i+4)})S(g_0^{(1)} g_1^{(2)} 
\dots g_n^{(n+1)})\cdot a_0 ,\\
 & S(g_1^{(1)} g_2^{(2)}\dots g_n^{(n)})\cdot a_1, 
 \dots,S(g_{i+1}^{(1)} \dots  g_n^{(n+1-i)})(S(g_i^{(1)})(a_i)a_{i+1}), \dots, \\
& S(g_{n-1}^{(1)} g_n^{(2)}) \cdot a_{n-1}, S(g_n^{(1)}) \cdot a_n ))=
\Phi d_i^{A\#_\sigma\mathcal{H}}( a_0 \otimes g_0, \dots , a_n \otimes g_n). 
\end{align*}
For $i=n$, we have:
 \begin{align*}
&  d_n \delta_n\Phi (a_0 \otimes g_0,\dots , a_n \otimes g_n)\\
&= d_n \delta_n(g_0^{(2)} , g_1^{(3)} , \dots ,g_n^{(n+2)} \mid S(g_0^{(1)} g_1^{(2)} \dots g_n^{(n+1)}) 
\cdot a_0, S(g_1^{(1)} g_2^{(2)} \dots g_n^{(n)}) \cdot a_1,\\
&  \hspace{7cm} \dots,
S(g_{n-1}^{(1)} g_n^{(2)}) \cdot a_{n-1},S(g_n^{(1)}) \cdot a_n)\\
&= d_n (g_0^{(3)} , g_1^{(4)} , \dots ,g_n^{(n+3)} \mid \\
&  \hspace{2cm}(S(g_0^{(2)} g_1^{(3)} \dots g_n^{(n+2)})S(g_n^{(1)}) \cdot a_n)(S(g_0^{(1)} g_1^{(2)} 
\dots g_n^{(n+1)}) \cdot a_0),\\
&  \hspace{5cm} S(g_1^{(1)} g_2^{(2)} \dots g_n^{(n)}) \cdot a_1, \dots,
S(g_{n-1}^{(1)} g_n^{(2)}) \cdot a_{n-1})\\
&= (g_n^{(2n+4)}g_0^{(3)} , g_1^{(4)} , \dots ,g_{n-1}^{(n+2)} \mid \\
&  \hspace{2cm} \sigma(g_n^{(2n+5)},g_0^{(4)})g_n^{(n+3)}\cdot((S(g_0^{(2)} g_1^{(3)} \dots g_n^{(n+2)})S(g_n^{(1)}) 
\cdot a_n)\\
&(S(g_0^{(1)} g_1^{(2)} \dots g_n^{(n+1)}) \cdot a_0)), S(g_1^{(1)} g_2^{(2)} \dots g_n^{(n)}) \cdot a_1, \dots,
S(g_{n-1}^{(1)} g_n^{(2)}) \cdot a_{n-1})\\
&= \Phi d_n^{A \#_\sigma \mathcal{H}} (a_0 \otimes g_0 , \dots , a_n \otimes g_n).
\end{align*}

The commutativity of $\Phi$
and  the degeneracies  are easier to check and is left to the reader.  The theorem is proved.    
\end{proof}
  Let $H_\sigma=k\underset{\sigma}{\#}\mathcal{H}$ denote the crossed product 
  of $\mathcal{H}$ and $k$ where $\mathcal{H}$ acts on $k$ via the counit $\epsilon$. One can check 
  that the $q$-th row of the cylindrical module $A\underset{\sigma}{\natural}\mathcal{H}$ is
   the  standard Hochschild complex 
  of the algebra $\mathcal{H}_\sigma$ with coefficients in the bimodule $M_q=\mathcal{H}\ot A^{\ot(q+1)}$.   
  Here $H_\sigma$ acts on $M_q$ on the left and right by 
  \begin{align*}
  &h\cdot(g\ot a_0\ot\dots \ot a_q)=\sigma(h^{(q+3)},g^{(2)})h^{(q+2)}g^{(1)}\ot h^{(1)}a_0\ot\dots\ot h^{(q+1)}a_q\\
  &(g\ot a_0\ot\dots \ot a_q)\cdot h=\sigma(g^{(2)},h^{(2)})g^{(1)}h^{(1)}\ot a_0\ot\dots\ot a_q.
   \end{align*}
   
   For the proof of Theorem 3.4 we need an extension of Mac Lane's isomorphism, which relates  group homology to 
   Hochschild homology, to Hopf algebras. 
   
   We recall that the Hopf homology of a Hopf algebra $\mathcal{H}$ 
    with coefficients in a left $\mathcal{H}$-module $M$ is  the  homology of the following complex 
    $$M\overset{d_0}{\leftarrow} \h\ot M\overset{d_1}{\leftarrow}
     \h\ot\h\ot M\overset{d_2}{\leftarrow}\dots\h^{\ot n}\ot M 
     \overset{d_n}{\leftarrow}\h^{\ot(n+1)}\ot M \leftarrow\dots,$$
     where the differential $d_n$ is given by 
     \begin{align*}
     &d_n(h_0\ot h_1\ot \dots\ot h_n\ot m)=
      \epsilon(h_0)h_1\ot\dots \ot h_n\ot m+\\&
      \sum_{1\le i\le n-1}(-1)^{i}h_0\ot \dots\ot h_ih_{i+1}\ot
       \dots\ot h_n\ot m +(-1)^nh_0\ot h_1\ot\dots \ot h_{n-1}\ot hm.
     \end{align*}
      We denote the $n$th  Hopf homology group of $\h$ with coefficients in $M$ by $H_n(\h ;M)$. 
  
  Let $M$ be an $\mathcal{H}_\sigma$-bimodule. We can  convert  $M$
 to a new left $\mathcal{H}$-module, $\widetilde{M}=M$, where the action of $\mathcal{H}$  on $\widetilde{M}$ 
 is  defined by 
   \begin{equation}\label{act}
   h\blacktriangleright m=\sigma^{-1}(S(h^{(2)}),h^{(3)})\overline{h^{(4)}}m\overline{S(h^{(1)})},
   \end{equation}
   where $\bar{h}$ denotes the image of $h$ in $\mathcal{H}_\sigma$ under the map $h\rightarrow 1\#h$. 
   Note that in the proof of the following lemma the cocommutativity of $\mathcal{H}$ is used. 
   \begin{lemma}
   Let $M$ be  an $\mathcal{H}_\sigma$-bimodule. Then by the above definition $M$ is a left 
   $\mathcal{H}$-module, i.e.,  $g\blacktriangleright(h\blacktriangleright m)=(gh)\blacktriangleright m$, 
    for all $g,h\in \mathcal{H}$ and $m\in M$. 
    \end{lemma}
   \begin{proof} We have 
   \begin{align*} 
    &g\blacktriangleright(h\blacktriangleright m)=
    g\bt(\sigma^{-1}(S(h^{(2)}),h^{(3)})\overline{h^{(4)}}m\overline{S(h^{(1)})})\\
   & =\sigma^{-1}(S(g\sw{2}),g\sw{3})\sigma^{-1}(S(h^{(2)}),h^{(3)})\overline{g\sw{4}}(\overline{h^{(4)}}
   m\overline{S(h^{(1)})})\overline{S(g\sw{1})}\\
   &=\sigma^{-1}(S(g\sw{3}),g\sw{4})\sigma^{-1}(S(h^{(3)}),h^{(4)})\sigma(g\sw{5},h\sw{5})\\
   &\sigma(S(h\sw{2}),S(g\sw{2}))\overline{g\sw{5}h\sw{5}}
   m\overline{S(g\sw{1}h^{(1)})}\\
   &=\sigma^{-1}(S(g\sw{3}),g\sw{4})\sigma^{-1}(S(h^{(3)}),h^{(4)})\sigma(g\sw{5},h\sw{5})
   \sigma(S(h\sw{2}),S(g\sw{2}))\\
   &\sigma(S(h\sw{9})S(g\sw{9}),g\sw{6}h\sw{6})
   \sigma^{-1}(S(h\sw{8}S(g\sw{8}),g\sw{7}h\sw{7})\overline{g\sw{5}h\sw{5}}
   m\overline{S(g\sw{1}h^{(1)})}\\
   &=\sigma^{-1}(S(g\sw{3}),g\sw{4})\sigma^{-1}(S(h^{(3)}),h^{(4)})\sigma(g\sw{5},h\sw{5})
   \sigma(S(h\sw{2}),S(g\sw{2}))\\
   &\sigma(S(h\sw{8})S(g\sw{8}),g\sw{6}h\sw{6})
   \sigma^{-1}(S(h\sw{9}S(g\sw{9}),g\sw{7}h\sw{7})\overline{g\sw{5}h\sw{5}}
   m\overline{S(g\sw{1}h^{(1)})}\\
   &=\sigma(S(h\sw{9})S(g\sw{9}),g\sw{5})\sigma(S(h\sw{8}),h\sw{5})\sigma^{-1}(S(g\sw{3}),g\sw{4})\sigma^{-1}(S(h^{(3)}),h^{(4)})
   \sigma(S(h\sw{2}),S(g\sw{2}))\\
   &\sigma^{-1}(S(h\sw{9}S(g\sw{9}),g\sw{7}h\sw{7})\overline{g\sw{5}h\sw{5}}
   m\overline{S(g\sw{1}h^{(1)})}\\
      &=\sigma(S(h\sw{5})S(g\sw{5}),g\sw{8})\sigma^{-1}(S(g\sw{5}),g\sw{6})
   \sigma(S(h\sw{6}),S(g\sw{6}))\\
   &\sigma^{-1}(S(h\sw{4}S(g\sw{4}),g\sw{3}h\sw{3})\overline{g\sw{2}h\sw{2}}
   m\overline{S(g\sw{1}h^{(1)})}\\
   &=\sigma^{-1}(S(g\sw{8}),g\sw{7})\sigma(S(g\sw{6}),g\sw{5})\sigma^{-1}(S(h\sw{4}S(g\sw{4}),g\sw{3}h\sw{3})
   \overline{g\sw{2}h\sw{2}}
   m\overline{S(g\sw{1}h^{(1)})}\\
      &=\sigma^{-1}(S(h\sw{2})S(g\sw{2}),g\sw{3}h\sw{3})
   \overline{g\sw{4}h\sw{4}}
   m\overline{S(g\sw{1}h^{(1)})}=(gh)\bt m.\\
      \end{align*}
   \end{proof}
The following result was first proved in \cite{kr2} for $\sigma$ a trivial cocycle. 
  \begin{theorem}[ Mac Lane Isomorphism for Hopf crossed products]
  Let $M$ be an  $\mathcal{H}_\sigma$-bimodule and $\widetilde{M} $ be defined as above. 
  Then the following map defines  an isomorphism  between Hochschild and Hopf homology complexes:
  \begin{center}
    $\Theta: C_n(\mathcal{H}_\sigma,M)\longrightarrow C_n(\mathcal{H};\widetilde{M})$ \\
    $\Theta(\bar{h}_1\ot\dots\ot \bar{h}_n\ot m)=h_1^{(2)}\ot h_2^{(2)}\ot\dots\ot h_n^{(2)}
    \ot m \overline{h_1^{(1)}}\dots\overline{h_n^{(1)}}.$ 
    \end{center}
  \end{theorem}
  \begin{proof}
  We show more than what we need for the proof,  namely  we show that $\Theta$ is an isomorphisms of 
   of simplicial modules. We have: 
   \begin{align*}
   &\Theta \delta_0(\bar{h}_1\ot\dots\ot \bar{h}_n\ot m)=\Theta(\bar{h}_2\ot\dots\ot \bar{h}_n\ot m\bar{h}_1)=\\
   &h_2^{(2)}\ot\dots\ot h_n^{(2)}
    \ot m \overline{h_1}\overline{h_2^{(1)}}\dots\overline{h_n^{(1)}}
     =\delta_0\Theta(\bar{h}_1\ot\dots\ot \bar{h}_n\ot m).  
   \end{align*}
    For $0\le i\le n$, we have 
   \begin{align*}
   &\Theta\delta_i(\bar{h}_1\ot\dots\ot \bar{h}_n\ot m)=\\
   &\Theta (\bar{h}_1\ot\dots \ot (\ov{h_i})( \ov{h_{i+1}})\ot \dots \ot \ov{h}_n\ot m)=\\
   &h_1^{(2)}\ot h_2^{(2)}\ot\dots\ot h_i^{(2)}h_{i+1}^{(2)} \ot\dots\ot h_n^{(2)}
    \ot m \overline{h_1^{(1)}}\dots\overline{h_n^{(1)}}=\\
    &\delta_i\Theta (\bar{h}_1\ot\dots\ot \bar{h}_n\ot m).  
 \end{align*}
   We leave it  to the reader to check the commutativity of $\Theta$ with the last face and the degeneracies.
   
   To finish the proof one can check that  the following map is the inverse of $\Theta$ 
   \begin{center}
    $\mathfrak{T}:C_n(\mathcal{H};\widetilde{M}) \longrightarrow C_n(\mathcal{H}_\sigma,M)$ \\
    $\mathfrak{T}(h_1\ot\dots\ot h_n\ot m)=\sigma^{-1}(S(h_1^{(2)}),h_1^{(3)})\dots \sigma^{-1}(S(h_n^{(2)}),h_n^{(3)})
    \ov{h_1^{(4)}}\ot\dots\ot\ov{h_n^{(4)}}\ot m\ov{S(h_n^{(1)})}\dots \ov{S(h_1^{(1)})}.$ 
    \end{center}
     \end{proof}
   
We  apply the generalized cyclic Eilenberg-Zilber theorem (Theorem 2.1) to the cylindrical 
module $A\natural_\sigma\mathcal{H}$ to derive a spectral sequence for
 the cyclic homology of $A\#_\sigma\mathcal{H}$.   We have 
 $$Tot(A\natural_\sigma\mathcal{H})\rightarrow d(A\natural_\sigma\mathcal{H})\cong(A\#_\sigma\mathcal{H})^\natural,$$
  where the first map is a quasi-isomorphism of mixed complexes given in Theorem 2.1 
  and the second map is the isomorphism given in Theorem 3.2. 
  We filter the mixed complex $Tot(A\natural_\sigma\mathcal{H})$ by sub mixed complexes
  $$F^i(Tot(A\natural_\sigma\mathcal{H}))_n=\bigoplus_{\underset{q\le i}{p+q=n}}A^{\ot(p+1)}\ot\mathcal{H}^{\ot(q+1)}.$$
  This gives us an spectral sequence that converges  to $HC_\bullet(A\#_\sigma\mathcal{H}).$ We can then apply 
  Theorem 3.3 to identify the $E^1$-term of this spectral sequence, i.e. the homology of rows, as Hopf 
  homologies of $\mathcal{H}$, with coefficients in $M_q=\mathcal{H}\ot A^{\ot(q+1)}$, where $\mathcal{H}$ acts 
  on $M_q$ by
  \begin{align}\label{action}
  &h\blacktriangleright(g\ot a_0\ot a_1\ot\dots\ot a_q)= 
    \sigma^{-1}(S(h^{(3)}),h^{(4)})\sigma(h^{(q+6)},g^{(1)})
   \sigma(h^{(q+7)}g^{(2)},S(h^{2}))\notag \\
  & h^{(q+8)}g^{(3)}S(h^{(1)})\ot
    h^{(5)}(a_0)\ot\dots\ot h^{(q+5)}(a_q). 
    \end{align}      
       This proves the following theorem. 
   \begin{theorem}\label{spec}
   There is a spectral sequence that converges to $HC_{p+q}(A\#_\sigma\mathcal{H})$.
    The $E^1$-term of this spectral   sequence is given by 
	$$E^1_{p,q}=H_p(\mathcal{H};M_q).$$
  \end{theorem}

   Given any cylindrical module $X=\{X_{p,q}\}_{p,q\ge 0}$, if we compute the Hochschild  homologies of  rows of 
   $X$ we obtain a new bigraded $k$-module $X'=\{X'_{p,q}\}_{p,q\ge 0}$. 
   We claim  that the columns of $X'$, i.e. $\{X'_{p,q}\}_{q\ge 0}$ form  a cyclic module for each $p\ge 0.$ 
   For some special cases one  can find the  proof in \cite{gj,ak}. The same proof, however, 
    works  in the general case.
    This observation  proves the following proposition. 
 \begin{proposition}
 The $p^{\text{th}}$column of $E^1$,  i.e. $\{ H_p(\mathcal{H};M_q)\}_{q\ge 0}$ is a cyclic module for each $p\ge 0$.
    \end{proposition}
    We denote the $p^{\text{th}}$ column of $E^1$ by $N_p$. One can observe that the induced differential $d^1$ on $E^1$
   is simply the differential  $b+B$ associated to the cyclic modules $N_p$. This finishes the proof of the following 
   theorem. 
  \begin{theorem} \label{4.6}The $E^2$ term of the spectral  sequence  in Theorem \ref{spec} is 
  $$E^2_{p,q}=HC_q(N_p).$$
  \end{theorem}
 
 Recall that if $\mathcal{H}$ is a semisimple Hopf algebra, then
 for any $\mathcal{H}$-module $M$, $H_i(\mathcal{H},M)=0$
  for $i\ge 1$. From this and Theorem \ref{spec},  we obtain the following corollary. 
  \begin{corollary}\label{cor}
  Let $\mathcal{H}$ be semisimple. Then the above spectral sequence collapses and we obtain
  $$HC_q(A\underset{\sigma}{\#}\mathcal{H})\cong HC_q(N_0).$$ 
    \end{corollary}
  Since $N_0=H_0(\mathcal{H},M_\bullet)$, we obtain 
  \begin{equation*}
  N_{0,q}=M_q^\mathcal{H}=(\mathcal{H}\ot A^{\ot(q+1)})^\mathcal{H},
  \end{equation*}
  where the action of $\mathcal{H}$ is defined by \ref{action}.

\end{document}